\documentclass{svmult}
\usepackage{amsmath,amssymb}
\def\Frac#1#2{\frac
{
 {\raise.6ex
 \hbox{$\displaystyle#1$}}
}
{
 {\lower.6ex
 \hbox{$\displaystyle#2$}}
 }
}
\usepackage{graphicx}        
\def\eqref#1{(\ref{#1})}

\def\dsp#1{\displaystyle{#1}}
\mainmatter  
\begin{document}
\title*{On the computation of moments of the partial non-central chi-squared distribution function}
\titlerunning{Moments of the partial non-central chi-squared distribution function}
\author{A. Gil\inst{1}, J. Segura\inst{2} \and N.M. Temme\inst{3}}
\institute{
	Departamento de Matem\'atica Aplicada y Ciencias de la Computaci\'on,
	ETSI Caminos, Canales y Puertos, Universidad de Cantabria, 39005-Santander, Spain.
        e-mail: amparo.gil@unican.es 
	\and
	Departamento de Matem\'aticas, Estad\'{\i}stica y Computaci\'on,
        Universidad de Cantabria, 39005 Santander, Spain.
        e-mail: javier.segura@unican.es 
        \and
        CWI, Science Park 123, 1098 XG Amsterdam, The Netherlands. e-mail: nico.temme@cwi.nl}
\maketitle
\begin{abstract}
Properties satisfied by the moments of the partial non-central chi-square
distribution function, also known as Nuttall Q-functions, and methods for computing
these moments are discussed in this paper. The Nuttall Q-function is involved
in the study of a variety of problems in different fields, as for example digital communications. 
\end{abstract}

\section{Introduction} \label{intro}

The non-central chi-square distribution function of probability appears in many applications. For example,
in radar communications it appears when computing the detection of signals in noise
using a square-law detector. Its cumulative distribution function is also
known as the generalized Marcum $Q-$ function, which is defined by using the integral representation
\begin{equation}
\label{eq:intro01}
Q_{\mu} (x,y)=\displaystyle x^{\frac12 (1-\mu)} \int_y^{+\infty} t^{\frac12 (\mu -1)} e^{-t-x} I_{\mu -1} \left(2\sqrt{xt}\right) \,dt,
\end{equation}
where $\mu >0$ and $I_\mu(z)$ is the modified Bessel function.

In radar problems, if the signal-to-noise power ratio is $x$ for the sum of $\mu$ independent 
samples of the output of a square-law detector, this integral gives the probability of
that the sum will be $y$ or more.
 
The complementary function of the generalized Marcum $Q-$ function is given by

\begin{equation}
\label{eq:intro02}
P_{\mu} (x,y)=\displaystyle x^{\frac12 (1-\mu)}\int_0^{y} t^{\frac12 (\mu -1)} e^{-t-x} I_{\mu -1} \left(2\sqrt{xt}\right) \,dt,
\end{equation}
and the following relation holds
\begin{equation}\label{eq:intro03}
P_{\mu}(x,y)+Q_{\mu} (x,y)=1.
\end{equation}

Methods and an algorithm for computing the functions $P_{\mu}(x,y)$ and $Q_{\mu}(x,y)$ are described
in \cite{Gil:2012:CMQ}. 

The $\eta$th moment of the partial non-central chi-square distribution function is given by

\begin{equation}
\label{intro04}
Q_{\eta,\mu} (x,y)=\displaystyle x^{\frac12 (1-\mu)} \int_y^{+\infty} t^{\eta+\frac12 (\mu -1)} e^{-t-x} I_{\mu -1} \left(2\sqrt{xt}\right) \,dt,
\end{equation}

In this manuscript, we give properties satisfied by the moments of the partial non-central chi-square
distribution functions and discuss methods for computing
these moments, also known as Nuttall Q-functions \cite{Nuttall:1972:AHN}.
There are several applications where these functions are involved as for example, the analysis of the outage probability
of wireless communication systems with a minimum signal power constraint \cite{Simon:2002:MKS}, to mention
just one example within the telecommunications field.

\section{Properties}

The Maclaurin series for the modified Bessel function reads 

\begin{equation}\label{eq:defImu}
I_\mu(z)=\left(\tfrac12z\right)^\mu\sum_{n=0}^\infty \frac{\left(\frac14z^2\right)^n}{n!\,\Gamma(\mu+n+1)}.
\end{equation}

By substituting this expression in the integral representation, we obtain the series expansion for the
$\eta$th moment of the non-central chi-square distribution function:

\begin{equation}\label{qseries}
Q_{\eta,\,\mu} (x,y)=e^{-x}\sum_{n=0}^{\infty} \Frac{x^n}{n!} \Frac{\Gamma(\eta +\mu +n,\,y)}{\Gamma (\mu +n)}
\end{equation}

This expansion is given in terms of one of the standard incomplete gamma functions  defined by

\begin{equation}\label{eq:igfs}
\Gamma (\mu,x)=\int_{x}^{+\infty} t^{\mu-1} e^{-t} \,dt.
\end{equation}

Introducing the factor $\Gamma(\eta +\mu +n)$ in \eqref{qseries}, the expansion
can be also given in terms of the  incomplete gamma function ratio $Q_{\mu} (y)$, defined by

\begin{equation}
\label{eq:defincgam}
Q_{\mu} (x)=\frac{\Gamma (\mu,x)}{\Gamma(\mu)},
\end{equation}
and for which algorithms are given in \cite{gil:2012:IGR}.

 The expansion for the
$\eta$th moment of the non-central chi-square distribution function in terms of incomplete gamma function ratios is given by

\begin{equation}\label{qseries2}
Q_{\eta,\,\mu} (x,y)=e^{-x}\sum_{n=0}^{\infty} \Frac{x^n}{n!}\Frac{\Gamma(\eta+\mu+n)}{\Gamma(\mu+n)} 
Q_{\eta +\mu +n}(y)
\end{equation}

The series representation can be computed by using the algorithms for the incomplete gamma ratios 
described in  \cite{gil:2012:IGR}. The recurrence relation

\begin{equation}\label{eq:serexp01}
Q_{\eta+\mu+1}(y)=Q_{\eta+\mu}(y)+\frac{y^{\eta+\mu} e^{-y}}{\Gamma(\eta+\mu+1)},
\end{equation}
is stable for $Q_\mu(y)$ in the forward direction, so the evaluation of the terms in the series for this function in 
(\ref{qseries2}) is rather easy.

A recurrence relation for the moments of the non-central chi-squared distribution function can be obtained considering
integration by parts in the integral in \eqref{intro04}, together with the relation 
$z^\mu I_{\mu-1}(z)=\frac{d}{dz}\left(z^\mu I_\mu(z)\right)$. This gives

\begin{equation}
\label{RecQ1}
\dsp{Q_{\eta,\mu}(x,y)=Q_{\eta,\mu+1} (x,y)-\eta Q_{\eta-1,\mu+1} -\left(\Frac{y}{x}\right)^{\mu/2}y^{\eta} 
e^{-x-y} I_{\mu} (2\sqrt{xy})}.
\end{equation}

When $\eta=0$, this recurrence reduces to a first order difference equation for the Marcum-Q function
(see, for instance,  \cite{Temme:1993:ANA} \footnote{We note that a factor $e^{-y}$ is missing in \cite[Eq.~(1.4)]{Temme:1993:ANA}.}).
The recurrence relation given in \eqref{RecQ1} can be used for testing, and it can be also used for computation, as we
describe later. 

\section{Computing moments using the series expansion}

 The series expansion given in \eqref{qseries} has been tested by using the recurrence relation of \eqref{RecQ1}
written in the form

\begin{equation}
\Frac{Q_{\eta,\mu+1}(x,y)}{\dsp{Q_{\eta,\mu}(x,y)+\eta Q_{\eta-1,\mu+1} +\left(\Frac{y}{x}\right)^{\mu/2}y^{\eta} 
e^{-x-y} I_{\mu} (2\sqrt{xy})} }=1.      
\label{errRR}
\end{equation}

The deviations from $1$ of the left-hand side of \eqref{errRR} (in absolute value) will measure the accuracy of the tested methods. 
The series expansion has been implemented in the Fortran 90 module {\bf NuttallF}. This module uses 
another module ({\bf IncgamFI}) for the computation of the gamma function ratios. 
We have tested the parameter region $(\eta,\,\mu,\,x,\,y) \in (1,\,50) \times (1,\,50) \times (0,\,20) \times (0,\,20)$.
The tests show that an accuracy better than $10^{-12}$ in this region can be obtained 
with the series expansion. 

When $\mu$ or $\mu +n$ are large, it is convenient to use approximations for the ratio of gamma functions appearing 
in the expression, in order to avoid the appearance of overflow problems sooner than expected. In the case $\mu +n \rightarrow \infty$
we have:

\begin{equation}
\Frac{\Gamma(\eta+\mu+n)}{\Gamma(\mu+n)}\sim (\mu+n)^{\eta}.
\end{equation}

The following table shows some values of moments
of the chi-square distribution function computed with the series expansion and the
corresponding values obtained with the direct computation of the integral representation
using Maple with 50 digits (the results shown in the table correspond to the first 18 digits
obtained with these computations). 
The computation of the series expansion has been implemented in the double precision Fortran 90 module 
 {\bf NuttallF}.
As can be seen, an agreement of minimum 14-15 digits is obtained in all cases, which is consistent
with the expected accuracy of the double precision Fortran 90 module.   

In some cases, Maple fails to compute the integral and acceleration can be obtained by suitably truncating
the improper integral and changing the variable of integration. We notice that, as before commented, the
modified Bessel function is exponentially increasing for large arguments and then the integrand in 
(\ref{intro04}) can be estimated by $t^\gamma e^{-(\sqrt{t}-\sqrt{x})^2}$, $\gamma =\eta+(\mu-1)/2$
which is related to a Gaussian centered $t=x$. The maximum value of this function is attained at 
$t=(\sqrt{x}+\sqrt{x+4\gamma})^2/4$ and integrating around this value with a sufficiently wide interval
is enough. This truncated integral over finite interval $[a,b]$ can be then transformed with
a linear change to an integral in $[-1,1]$ and the convergence is further accelerated by considering
the change of variable $t=\tanh (u)$, particularly if the trapezoidal rule is used
for evaluating the integral (see \cite[\S5.4.2]{Gil:2007:NSF}). These modifications are observed to speed up
the computation of the integrals using Maple, particularly for the last value in Table 1 for which 
Maple does not appear to be able to converge to an accurate value.

\begin{table}
$$
\begin{array}{c|c|c|c|c|c}
\eta & \mu & x & y &  Q_{\eta,\mu}(x,y) & Q_{\eta,\mu}(x,y)\mbox{ with Maple} \\
\hline
1 & 1 & 0.1 & 1.5 &   0.6644091427683566 &          0.66440914276835656 \\
5 & 10 & 0.1 & 1.5 &  252472.22699183668   &        252472.226991836658  \\                 
50 & 30 & 0.1 & 1.5 &  1.1944632251434243\,10^{+86}    &  1.19446322514344860\,10^{+86}  \\
1 & 1 & 1.2 & 5 & 0.5457546041478581 &              0.54575460414785805       \\
5 & 10 & 1.2 & 5 &   419098.1927146542 &            419098.192714654143 \\
50 & 30 & 1.2 & 5 &  6.809314196073125\,10^{+86} &        6.80931419607285639 \,10^{+86} \\
1 & 1 & 5 & 10 & 1.4822515303982464  &              1.48225153039824667  \\
5 & 10 & 5 & 10 & 1654969.264263704  &              1654969.26426370245  \\
50 & 30 & 5 & 10 &  1.1734657613338925\,10^{+89}  & 1.17346576133388184 \,10^{+89}
 \\
\end{array}
$$
{\footnotesize {\bf Table 1.} Values of the moments
of the chi-square distribution function for different choices of the parameters
$\eta$, $\mu$ $x$ and $y$. The values shown 
are obtained with the series expansion and with the direct computation of the integral representation
using Maple with 50 digits.}
\label{table1}
\end{table}

\section{Computing moments by recursion}
If we write the recurrence relation (\ref{RecQ1}) as
\begin{equation}
\label{rec2}
Q_{\eta,\mu+1} (x,y)=\dsp{Q_{\eta,\mu}(x,y)+\eta Q_{\eta-1,\mu+1} +\left(\Frac{y}{x}\right)^{\mu/2}y^{\eta} 
e^{-x-y} I_{\mu} (2\sqrt{xy})}
\end{equation}
then it is clear that we have a numerically stable relation because all the terms in the right hand side are positive.

Now, assume that the moments of order zero (Marcum functions) $Q_{0,\mu}$ are known for $\mu=1,2,\ldots, N$ (or for a
sequence of real values $\mu_i$, $i=1,\ldots N$, with $\mu_{i+1}-\mu_i=1$). If $Q(1,\mu)$ is also known, the relation
(\ref{rec2}) can be used to compute $Q(1,\mu +1)$; therefore, starting from the value $Q(1,1)$ we can compute 
$Q(1,\mu)$, $\mu=1,2,\ldots, N$ in a stable way. In the same way, after determining $Q(1,\mu)$, $\mu=1,2,\ldots, N$ and if
$Q(2,1)$ is know, we can compute $Q(1,\mu)$, $\mu=1,2,\ldots N$ and so on.

It is worth mentioning that the inhomogeneous recurrence has to be applied with care, particularly the inhomogeneous term.
As $x$ and/or $y$ becomes large the Bessel function increases exponentially; therefore we have the product of a small exponential times an exponentially large function and because of the bad conditioning of the exponentials, this translates into larger
relative errors; additionally, the exponentials may overflow/underflow. 
Part of this error can be avoided by considering the scaled Bessel function $\tilde{I}_{\nu}(x)=e^{-x}I_{\nu}(x)$.
In terms of this function
\begin{equation}
\label{rec2}
Q_{\eta,\mu+1} (x,y)=\dsp{Q_{\eta,\mu}(x,y)+\eta Q_{\eta-1,\mu+1} +\left(\Frac{y}{x}\right)^{\mu/2}y^{\eta} 
e^{-(\sqrt{x}-\sqrt{y})^2} \tilde{I}_{\mu} (2\sqrt{xy})}.
\end{equation}  

An alternative way of computing with recurrences is considering a homogeneous equation, which we can be constructed from
the inhomogeneous equation writing
\begin{equation}
\begin{array}{ll}
Q_{\eta,\mu+2}- Q_{\eta,\mu+1}-\eta Q_{\eta-1,\mu+2}=&c_{\mu+1}(Q_{\eta,\mu+1}- Q_{\eta,\mu}-\eta Q_{\eta-1,\mu+1}),\\[8pt]
& c_{\mu+1}=\displaystyle\sqrt{\Frac{y}{x}}\Frac{I_{\mu +1}(2\sqrt{xy})}{I_{\mu}(2\sqrt{xy})}.
\end{array}
\end{equation}
Then, if $Q(\eta-1,\mu)$ is know $\mu=1,2,\ldots N$, we can compute $Q(\eta,\mu)$ $\mu=1,2,\ldots N$, starting from
$Q(\eta,1)$ and $Q(\eta,2)$ with the recurrence
\begin{equation}
\label{rec4}
Q_{\eta,\mu+2}= (1+c_{\mu +1}) Q_{\eta,\mu+1}-c_{\mu +1}Q_{\eta,\mu}+\eta Q_{\eta-1,\mu+2}-\eta c_{\mu +1} Q_{\eta-1,\mu+1}.
\end{equation}

The advantage of this recurrence is that the overflow problems are reduced because ratios of Bessel functions appear instead
of Bessel functions themselves. Also, for computing these ratios, continued fraction representations
can be used. In Table 2 the use of the recurrence relation for computing $Q_{2,N}$ is tested for several values of $N$.
The values of $x$ and $y$ are fixed to $2$ and $3$, respectively. The table 
shows the relative
error obtained when comparing the value obtained with the recurrence relation and the direct computation using the series
expansion of \eqref{qseries2}:

  \begin{equation}
   \label{er}
    E_r=\left|1-\Frac{Q^S_{2,N}(2,3)}{Q^R_{2,N}(2,3)}\right|.
   \end{equation}

The continued fraction for the ratio of Bessel functions is computed using the modified Lentz algorithm 
\cite{Thompson:1986:CFL} and  \cite[\S6.6.2]{Gil:2007:NSF}.

\begin{table}
$$
\begin{array}{c|c|}
N &  \mbox{Relative error } \eqref{er} \\
\hline
10 &  \, 2.96\,10^{-16}       \\
20 &  \,4.84\,10^{-16}   \\
30 &  \,1.67\,10^{-15}  \\
40 &  \,1.02\,10^{-15} \\
50 &  \,9.02\,10^{-15}   \\ 
60 &  \,3.09\,10^{-15}   \\
\end{array}
$$
{\footnotesize {\bf Table 2} Test of the application of the recurrence relation given in \eqref{rec4}. 
The relative errors are obtained when comparing the value obtained with the recurrence relation \eqref{rec4} and the direct computation using the series expansion of \eqref{qseries2}.}
\end{table}

\section*{Acknowledgements}
This work was supported by  {\emph{Ministerio de Ciencia e Innovaci\'on}}, 
project MTM2009-11686 and {\emph{Ministerio de Econom\'{\i}a y Competitividad}}, project
MTM2012-34787. 

\bibliographystyle{plain}

\end{document}